\documentclass[a4paper,11pt]{amsart}
\usepackage{graphicx}
\usepackage{mathptmx}
\usepackage{amsmath}
\usepackage{amssymb}
\usepackage{enumitem}
\usepackage{xcolor}

\newmuskip\pFqmuskip

\newcommand*\pFq[6][8]{%
  \begingroup 
  \pFqmuskip=#1mu\relax
  \mathcode`=\string"8000
  \begingroup\lccode`\~=`\,
  \lowercase{\endgroup\let~}\pFqcomma
  F^{#2}_{#3}{\left(\genfrac..{0pt}{}{#4}{#5}\bigg|#6\right)}%
  \endgroup
}
\newcommand{\pFqcomma}{\mskip\pFqmuskip}

\newtheorem{theorem}{Theorem}

\newtheorem{corollary}[theorem]{Corollary}

\begin{document}

\title[Jindalrae and Gaenari numbers and polynomials]{Jindalrae and Gaenari numbers and polynomials}

\author{Taekyun  Kim$^{1}$}
\address{$^{1}$
School of Science, Xi’an Technological University, Xi’an 710021, China\\
Department of Mathematics,  Kwangwoon University, Seoul 139-701, Republic of Korea}
\email{tkkim@kw.ac.kr}

\author{Dae San Kim$^{2}$}
\address{$^{2}$Department of Mathematics, Sogang University, Seoul 121-742, Republic of Korea}
\email{dskim@sogang.ac.kr}

\author{Lee-Chae Jang$^{3}$}
\address{$^{3}$Graduate School of Education, Konkuk University, Seoul 05029, Korea}
\email{lcjang@konkuk.ac.kr}

\author{Hyunseok  Lee$^{4}$}
\address{$^{4}$Department of Mathematics, Kwangwoon University, Seoul 139-701, Republic of Korea}
\email{luciasconstant@gmail.com}


\maketitle

\begin{abstract}
The aim of this paper is to study Jindalrae and Gaenari numbers and polynomials in connection with Jindalrae-Stirling numbers of the first and second kinds. For this purpose, we first introduce Jindalrae-Stirling numbers of the first and second kinds as extensions of the notions of the degenerate Stirling numbers of the first and second kinds, and deduce several relations involving those special numbers. Then we introduce Jindalrae and Gaenari numbers and polynomials and obtain some explicit expressions and identities associated with those numbers and polynomials. In addition, we interpret our results by using umbral calculus.

\end{abstract}

\section{Introduction and preliminaries}
Let $n$ be a nonnegaitve integer. Then the Stirling numbers of the first kind are defined as
\begin{equation}
	(x)_{n}=\sum_{l=0}^{n}S_{1}(n,l)x^{l},\quad (\mathrm{see}\ [5]), \label{1}
\end{equation}
where $(x)_{0}=1,\ (x)_{n}=x(x-1)\cdots(x-n+1)$, $(n\ge 1)$. \\
Whereas the Stirling numbers of the second kind are given by
\begin{equation}
	x^{n}=\sum_{l=0}^{n}S_{2}(n,l)(x)_{l},\quad (\mathrm{see}\ [10,19]). \label{2}
\end{equation}
By \eqref{1} and \eqref{2}, we get
\begin{equation}
\frac{1}{k!}\big(e^{t}-1\big)^{k}=\sum_{n=k}^{\infty}S_{2}(n,k)\frac{t^{n}}{n!},\quad (k\ge 0),\label{3}
\end{equation}
and
\begin{equation}
\frac{1}{k!}\big(\log(1+t)\big)^{k}=\sum_{n=k}^{\infty}S_{1}(n,k)\frac{t^{n}}{n!},\quad(k\ge 0).\label{4}
\end{equation}
\indent It is well known that the Bell polynomials are defined as
\begin{equation}
	e^{x(e^{t}-1)}=\sum_{n=0}^{\infty}B_{n}(x)\frac{t^{n}}{n!},\quad(\mathrm{see}\ [13,18]).\label{5}
\end{equation}
When $x=1,\ B_{n}=B_{n}(1)$ are called the Bell numbers. \\
\indent For $0 \neq \lambda\in\mathbb{R}$, the degenerate exponential functions are defined by
\begin{equation}
e_{\lambda}^{x}(t)=\sum_{n=0}^{\infty}\frac{(x)_{n,\lambda}}{n!}t^{n},\quad e_{\lambda}(t)=e_{\lambda}^{1}(t)=\sum_{n=0}^{\infty}\frac{(1)_{n,\lambda}}{n!}t^{n},\quad(\mathrm{see}\ [12,21]), \label{6}
\end{equation}
where $(x)_{0,\lambda}=1$, $(x)_{n,\lambda}=x(x-\lambda)(x-2\lambda)\cdots\big(x-(n-1)\lambda\big)$, $(n\ge 1)$. \\
\indent Let $\log_{\lambda}(t)$ be the compositional inverse of $e_{\lambda}(t)$, called the degenerate logarithm function, such that $\log_{\lambda}\big(e_{\lambda}(t)\big)=e_{\lambda}\big(\log_{\lambda}t\big)=t$. \\
Then we note that
\begin{equation}
\log_{\lambda}(1+t)=\frac{1}{\lambda}((1+t)^{\lambda}-1)=\sum_{n=1}^{\infty}\lambda^{n-1}(1)_{n,\frac{1}{\lambda}}\frac{t^{n}}{n!},\quad(\mathrm{see}\ [11]). \label{7}
\end{equation}
By \eqref{7}, we get $\displaystyle\lim_{\lambda\rightarrow 0}\log_{\lambda}(1+t)=\log(1+t)\displaystyle$. \\
\indent In [11], the degenerate Stirling numbers of the first kind are defined by
\begin{equation}
	(x)_{n}=\sum_{l=0}^{n}S_{1,\lambda}(n,l)(x)_{l,\lambda}. \label{8}
\end{equation}
As an inversion formula of \eqref{8}, the degenerate Stirling numbers of the second kind are defined by
\begin{equation}
	(x)_{n,\lambda}=\sum_{k=0}^{n}S_{2,\lambda}(n,k)(x)_{k},\ (n\ge 0),\ (\mathrm{see}\ [16]). \label{9}
\end{equation}
By \eqref{8} and \eqref{9}, we get
\begin{equation}
\frac{1}{k!}\big(\log_{\lambda}(1+t)\big)^{k}=\sum_{n=k}^{\infty}S_{1,\lambda}(n,k)\frac{t^{n}}{n!},\ (k\ge 0),\quad(\mathrm{see}\ [11]).\label{10}
\end{equation}
and
\begin{equation}
\frac{1}{k!}\big(e_{\lambda}(t)-1\big)^{k}=\sum_{n=k}^{\infty}S_{2,\lambda}(n,k)\frac{t^{n}}{n!},\ (k\ge 0),\quad(\mathrm{see}\ [16]),\label{11}
\end{equation}
\indent We define the degenerate Bell polynomials $B_{n,\lambda}(x)$ by
\begin{equation}
e_{\lambda}^{x}\big(e_{\lambda}(t)-1\big)=\sum_{n=0}^{\infty}B_{n,\lambda}(x)\frac{t^{n}}{n!}. \label{12}
\end{equation}
When $x=1$, $B_{n,\lambda}=B_{n,\lambda}(1)$ are called the degenerate Bell numbers. \\
From \eqref{6} and \eqref{10}, we note that
\begin{align}
	e_{\lambda}^{x}\big(e_{\lambda}(t)-1\big)\ &=\ \sum_{k=0}^{\infty}(x)_{k,\lambda}\frac{1}{k!}\big(e_{\lambda}(t)-1\big)^{k} \label{13} \\
	&=\  \sum_{k=0}^{\infty}(x)_{k,\lambda}\sum_{n=k}^{\infty}S_{2,\lambda}(n,k)\frac{t^{n}}{n!} \nonumber \\
	&=\ \sum_{n=0}^{\infty}\bigg(\sum_{k=0}^{n}(x)_{k,\lambda}S_{2,\lambda}(n,k)\bigg)\frac{t^{n}}{n!} \nonumber
\end{align}
By \eqref{12} and \eqref{13}, we get
\begin{equation}
B_{n,\lambda}(x)=\sum_{k=0}^{n}S_{2,\lambda}(n,k)(x)_{k,\lambda},\ (n\ge 0). \label{14}
\end{equation}
\indent Here we note that the so-called new type degenerate Bell polynomials $Bel_{n,\lambda}(x)$, which are different from the degenerate Bell polynomials just introduced, were considered recently in [18]. They are defined by the generating function $e_{\lambda}^{x}\big(e^t -1\big)=\sum_{n=0}^{\infty}Bel_{n,\lambda}(x)\frac{t^{n}}{n!}$, so that $Bel_{n,\lambda}(x)=\sum_{k=0}^{n}S_2(n,k)(x)_{k,\lambda},\ (n\ge 0)$.

\vspace{0.1in}

In [2,3], Carlitz initiated a study of the degenerate Bernoulli and Euler polynomials and numbers, which are degenerate versions of the usual Bernoulli and Euler polynomials and numbers (see [1,22,25,30,31]). In recent years, studying degenerate versions of quite a few special polynomials and numbers regained lively interest of some mathematicians and yielded many interesting results (see [4-7,10-13,15-18,20,21,23,24,26,27,29]). Here we note that such degenerate versions of many special polynomials and numbers have been investigated by employing different tools like generating functions, combinatorial methods, umbral calculus techniques, probability theory, $p$-adic analysis and differential equations, etc.\\
\indent The aim of the present paper is to study Jindalrae and Gaenari polynomials and
numbers in connection with Jindalrae-Stirling numbers of the first and second kinds, and find some arithmetic and combinatorial results on those polynomials and numbers. First, we define Jindalrae-Stirling numbers of the first and second kinds, as extensions of the notions of the degenerate Stirling numbers of the first and second kinds, and find some relations involving those special numbers. Then we introduce Jindalrae and Gaenari numbers and polynomials and obtain some explicit expressions and identities associated with those numbers and polynomials. In addition, we interpret our results by using umbral calculus.\\
\indent This paper is organized as follows. In Section 1, we go over some necessary stuffs that are needed throughout this paper. These include the degenerate exponential functions, the degenerate logarithm function, the degenerate Stirling numbers of the first and second kinds and the degenerate Bell numbers. Here we note that
the degenerate Bell polynomials $B_{n,\lambda}(x)$ (see \eqref{12}) are different from the partially degenerate Bell polynomials $bel_{n,\lambda}(x)$ in [15] and also from the new type degenerate Bell polynomials $Bel_{n,\lambda}(x)$ in [18]. In Section 2, we introduce Jindalrae-Stirling numbers of the first and second kinds, as extensions of the notions of the degenerate Stirling numbers of the first and second kinds and find some relations connecting those special numbers, the degenerate Stirling numbers of the first and second kinds and the degenerate Bell numbers and polynomials. Then we introduce Jindalrae numbers and polynomials, as an extension of the notion of the degenerate Bell numbers and polynomials, and Gaenari numbers and polynomials, and find some explicit expressions and identities involving those numbers and polynomials, Jindalrae-Stirling numbers of the first and second kinds, the degenerate Stirling numbers of the first and second kinds and the degenerate Bell polynomials.  In Section 3, we interpret the results in Section 2 by means of umbral calculus. Finally, we conclude this paper in Section 4. \\

\vspace{0.1in}

\section{Jindalrae and Gaenari numbers and polynomials}
By replacing $t$ by $e^{t}-1$ in \eqref{3}, we get
\begin{align}
	\frac{1}{k!}\big(e^{e^{t}-1}-1\big)^{k}\ &=\ \sum_{m=k}^{\infty}S_{2}(m,k)\frac{1}{m!}(e^{t}-1)^{m} \label{15} \\
	&=\ \sum_{m=k}^{\infty}S_{2}(m,k)\sum_{n=m}^{\infty}S_{2}(n,m)\frac{t^{n}}{n!} \nonumber \\
	&=\ \sum_{n=k}^{\infty}\bigg(\sum_{m=k}^{n}S_{2}(m,k)S_{2}(n,m)\bigg)\frac{t^{n}}{n!} \nonumber
\end{align}
Let
\begin{equation} \frac{1}{k!}\big(e^{e^{t}-1}-1\big)^{k}=\sum_{n=k}^{\infty}T(n,k)\frac{t^{n}}{n!}. \label{16}
\end{equation}
Then, by \eqref{15} and \eqref{16}, we get
\begin{equation}
	T(n,k)=\sum_{m=k}^{\infty}S_{2}(n,m)S_{2}(m,k), \label{17}
\end{equation}
where $n,k\ge 0$, with $n\ge k$. \\
\indent Further, we have
\begin{align}
	\frac{1}{k!}\big(e^{e^{t}-1}-1\big)^{k}\ &=\ \frac{1}{k!}\bigg(\sum_{n=1}^{\infty}B_{n}\frac{t^{n}}{n!}\bigg)^{k}\label{18} \\
	&=\ \frac{1}{k!}\sum_{n=k}^{\infty}\bigg(\sum_{n_{1}+\cdots+n_{k}=n}\binom{n}{n_{1},\dots,n_{k}}B_{n_{1}}\cdots B_{n_{k}}\bigg)\frac{t^{n}}{n!},  \nonumber
\end{align}
where the inner sum runs over all positive integers $n_1,\dots,n_k$, with $n_1+\cdots+n_k=n$.
Thus, by \eqref{16}, \eqref{17} and \eqref{18}, we get
\begin{align}\label{19}
	T(n,k)\ &=\ \frac{1}{k!} \sum_{n_{1}+\cdots+n_{k}=n}\binom{n}{n_{1},\dots,n_{k}}B_{n_{1}}\cdots B_{n_{k}}\\
	&=\ \sum_{m=k}^{n}S_{2}(n,m)S_{2}(m,k). \nonumber
\end{align}
Note that $\displaystyle T(n,1)=\sum_{m=1}^{n}S_{2}(n,m)=B_{n},\ (n\ge 1)$. \\
\indent For $k\ge 0$,  as an extension of the notion of the degenerate Stirling numbers of the second kind we define Jindalrae-Stirling numbers of the second kind (see \eqref{9} and \eqref{11}) by
\begin{equation}
	\frac{1}{k!}\big(e_{\lambda}(e_{\lambda}(t)-1)-1\big)^{k}=\sum_{n=k}^{\infty}S_{J,\lambda}^{(2)}(n,k)\frac{t^{n}}{n!}. \label{20}
\end{equation}
On the other hand,
\begin{align}
	\frac{1}{k!}\big(e_{\lambda}(e_{\lambda}(t)-1)-1\big)^{k}\ &=\ \sum_{m=k}^{\infty}S_{2,\lambda}(m,k)\sum_{n=m}^{\infty}S_{2,\lambda}(n,m)\frac{t^{n}}{n!}\label{21} \\
	&=\ \sum_{n=k}^{\infty}\bigg(\sum_{m=k}^{n}S_{2,\lambda}(m,k)S_{2,\lambda}(n,m)\bigg)\frac{t^{n}}{n!}. \nonumber
\end{align}
Therefore, by \eqref{20} and \eqref{21}, we obtain the following theorem.
\begin{theorem}
For $n,k\ge 0$, with $n\ge k$, we have
\begin{displaymath}
S_{J,\lambda}^{(2)}(n,k)=\sum_{m=k}^{n}S_{2,\lambda}(n,m)S_{2,\lambda}(m,k).
\end{displaymath}
\end{theorem}

When $k=1$, we have
\begin{align*}
	\sum_{n=1}^{\infty}S_{J,2}^{(2)}(n,1)\frac{t^{n}}{n!}\ &=\ e_{\lambda}\big(e_{\lambda}(t)-1\big)-1 \\
	&=\ \sum_{n=0}^{\infty}B_{n,\lambda}\frac{t^{n}}{n!}-1\ =\ \sum_{n=1}^{\infty}B_{n,\lambda}\frac{t^{n}}{n!}.
\end{align*}
For $n\ge 1$, we have
\begin{align}
	S_{J,\lambda}^{(2)}(n,1)\ =\ B_{n,\lambda}\ &=\ \sum_{m=1}^{n}S_{2,\lambda}(n,m)S_{2,\lambda}(m,1) \label{22} \\
	&=\ \sum_{m=1}^{n}S_{2,\lambda}(n,m)(1)_{m,\lambda}. \nonumber
\end{align}
For $k\ge 0$, by replacing $t$ by $\log_{\lambda}(1+t)$ in \eqref{20}, we get
\begin{align}
\frac{1}{k!}\big(e_{\lambda}(t)-1\big)^{k}\ &=\ \sum_{m=k}^{\infty}S_{J,\lambda}^{(2)}(m,k)\frac{1}{m!}\big(\log_{\lambda}(1+t)\big)^{m} \label{23} \\	
&=\ \sum_{m=k}^{\infty}S_{J,\lambda}^{(2)}(m,k)\sum_{n=m}^{\infty}S_{1,\lambda}(n,m)\frac{t^{n}}{n!} \nonumber \\
&=\ \sum_{n=k}^{\infty}\bigg(\sum_{m=k}^{n}S_{J,\lambda}^{(2)}(m,k)S_{1,\lambda}(n,m)\bigg)\frac{t^{n}}{n!}. \nonumber
\end{align}
Therefore, by \eqref{10} and \eqref{23}, we obtain the following theorem.
\begin{theorem}
	For $n,k\ge 0$, with $n\ge k$, we have
	\begin{displaymath}
		S_{2,\lambda}(n,k)=\sum_{m=k}^{n}S_{J,\lambda}^{(2)}(m,k)S_{1,\lambda}(n,m).
	\end{displaymath}
\end{theorem}
When $k=1$, we have
\begin{align}
	S_{2,\lambda}(n,1)\ &=\ \sum_{m=1}^{n}S_{J,\lambda}^{(2)}(m,1)S_{1,\lambda}(n,m) \label{24} \\
	&=\ \sum_{m=1}^{n}B_{m,\lambda}S_{1,\lambda}(n,m). \nonumber
\end{align}
It is easy to show that $S_{2,\lambda}(n,1)=(1)_{n,\lambda}$,\, $(n \geq 1)$. Therefore, by \eqref{24}, we obtain the following corollary.
\begin{corollary}
	For $n\in\mathbb{N}$, we have
	\begin{displaymath}
		\sum_{m=1}^{n}B_{m,\lambda}S_{1,\lambda}(n,m)=(1)_{n,\lambda}.
	\end{displaymath}
\end{corollary}
As an inversion formula of \eqref{20} and an extension of the notion of the degenerate Stirling numbers of the first kind (see \eqref{8} and \eqref{10}), we define Jindalrae-Stirling numbers of the first kind by
\begin{equation}
	\frac{1}{k!}\big(\log_{\lambda}(\log_{\lambda}(1+t)+1)\big)^{k}=\sum_{n=k}^{\infty}S_{J,\lambda}^{(1)}(n,k)\frac{t^{n}}{n!},\quad(k\ge 0). \label{25}
\end{equation}
We note that
\begin{align}
	\frac{1}{k!}\big(\log_{\lambda}(\log_{\lambda}(1+t)+1)\big)^{k}\ &=\ \sum_{m=k}^{\infty}S_{1,\lambda}(m,k)\frac{1}{m!}\big(\log_{\lambda}(1+t)\big)^{m} \label{26} \\
	&=\ \sum_{m=k}^{\infty}S_{1,\lambda}(m,k)\sum_{n=m}^{\infty}S_{1,\lambda}(n,m)\frac{t^{n}}{n!} \nonumber \\
	&=\ \sum_{n=k}^{\infty}\bigg(\sum_{m=k}^{n}S_{1,\lambda}(m,k)S_{1,\lambda}(n,m)\bigg)\frac{t^{n}}{n!}. \nonumber
\end{align}
Therefore, by \eqref{25} and \eqref{26}, we obtain the following theorem.
\begin{theorem}
	For $n,k\ge 0$, with $n\ge k$, we have
	\begin{displaymath}
		S_{J,\lambda}^{(1)}(n,k)=\sum_{m=k}^{n}S_{1,\lambda}(n,m)S_{1,\lambda}(m,k).
	\end{displaymath}
\end{theorem}
When $k=1$, we have
\begin{align}
S_{J,\lambda}^{(1)}(n,1)\ &=\ \sum_{m=1}^{n}S_{1,\lambda}(m,1)S_{1,\lambda}(n,m) \label{27} \\
&=\ \sum_{m=1}^{n}(m-1)!\binom{\lambda-1}{m-1}S_{1,\lambda}(n,m). \nonumber
\end{align}
\begin{corollary}
	For $n\ge 1$, we have
	\begin{displaymath}
	S_{J,\lambda}^{(1)}(n,1)=\sum_{m=1}^{n}(m-1)!\binom{\lambda-1}{m-1}S_{1,\lambda}(n,m).
	\end{displaymath}
\end{corollary}
From \eqref{20}, we note that
\begin{align}
	\sum_{n=k}^{\infty}S_{J,\lambda}^{(2)}(n,k)\frac{t^{n}}{n!}\ &=\ \frac{1}{k!}\big(e_{\lambda}(e_{\lambda}(t)-1)-1\big)^{k} \label{28} \\
	&=\ \frac{1}{k!}\sum_{l=0}^{k}\binom{k}{l}(-1)^{k-l}e_{\lambda}^{l}\big(e_{\lambda}(t)-1\big) \nonumber \\
	&=\ \sum_{n=0}^{\infty}\bigg(\frac{1}{k!}\sum_{l=0}^{k}\binom{k}{l}(-1)^{k-l}B_{n,\lambda}(l)\bigg)\frac{t^{n}}{n!}, \nonumber
\end{align}
where $k$ is a nonnegative integer. \\
By comparing the coefficients on both sides of \eqref{28}, we get
\begin{equation}
	\label{29} \frac{1}{k!}\sum_{l=0}^{k}\binom{k}{l}(-1)^{k-l}B_{n,\lambda}(l)=\left\{\begin{array}{ccc}
	S_{J,\lambda}^{(2)}(n,k), & \textrm{if $n\ge k$,}\\
	0, & \textrm{if $ 0 \leq n <k$.}
\end{array}\right.
\end{equation}
Therefore, by \eqref{29}, we obtain the following theorem.
\begin{theorem}
	For $n,k\ge 0$, with $n\ge k$, we have
	\begin{displaymath}
	S_{J,\lambda}^{(2)}(n,k)=\frac{1}{k!}\sum_{l=0}^{k}\binom{k}{l}(-1)^{k-l}B_{n,\lambda}(l).
	\end{displaymath}
\end{theorem}
Now, we observe that
\begin{align}
	\sum_{k=0}^{\infty}(x)_{k,\lambda}\frac{1}{k!}\big(\log_{\lambda}(\log_{\lambda}(1+t)+1)\big)^{k}\ &=\  \sum_{k=0}^{\infty}(x)_{k,\lambda}\sum_{n=k}^{\infty}S_{J,\lambda}^{(1)}(n,k)\frac{t^{n}}{n!} \label{30} \\
	&=\ \sum_{n=0}^{\infty}\bigg(\sum_{k=0}^{n}(x)_{k,\lambda}S_{J,\lambda}^{(1)}(n,k)\bigg)\frac{t^{n}}{n!}. \nonumber
\end{align}
On the other hand,
 \begin{equation}
 	\sum_{k=0}^{\infty}(x)_{k,\lambda}\frac{1}{k!}\big(\log_{\lambda}(\log_{\lambda}(1+t)+1)\big)^{k}\ =\ e_{\lambda}^{x}(\log_{\lambda}(\log_{\lambda}(1+t)+1)\big). \label{31}
\end{equation}
\begin{align*}
	&=\ \big(\log_{\lambda}(1+t)+1\big)^{x}\ =\ \sum_{l=0}^{\infty}(x)_{l}\frac{1}{l!}\big(\log_{\lambda}(1+t)\big)^{l}\\
	&=\ \sum_{l=0}^{\infty}(x)_{l}\sum_{n=l}^{\infty}S_{1,\lambda}(n,l)\frac{t^{n}}{n!}\ =\ \sum_{n=0}^{\infty}\bigg(\sum_{l=0}^{n}S_{1,\lambda}(n,l)(x)_{l}\bigg)\frac{t^{n}}{n!}.
\end{align*}
Thus, by \eqref{30} and \eqref{31}, we get
\begin{equation}
\sum_{l=0}^{n}S_{1,\lambda}(n,l)(x)_{l}=\sum_{k=0}^{n}S_{J,\lambda}^{(1)}(n,k)(x)_{k,\lambda},\ (n\ge 0). \label{32}
\end{equation}
From \eqref{9} and \eqref{32}, we can derive the following equation (33),
\begin{align}
	\sum_{l=0}^{n}S_{1,\lambda}(n,l)(x)_l\ &=\ \sum_{k=0}^{n}S_{J,\lambda}^{(1)}(n,k)(x)_{k,\lambda} \label{33}\\
	&=\ \sum_{k=0}^{n}S_{J,\lambda}^{(1)}(n,k)\sum_{l=0}^{k}S_{2,\lambda}(k,l)(x)_{l} \nonumber \\
	&=\ \sum_{l=0}^{n}\bigg(\sum_{k=l}^{n}S_{J,\lambda}^{(1)}(n,k)S_{2,\lambda}(k,l)\bigg)(x)_{l} \nonumber
\end{align}
Therefore, by comparing the coefficients as both sides of \eqref{33}, we obtain the following theorem.
\begin{theorem}
	For $n,l\ge 0$, we have
	\begin{displaymath}
		S_{1,\lambda}(n,l)=\sum_{k=l}^{n}S_{J,\lambda}^{(1)}(n,k)S_{2,\lambda}(k,l).
	\end{displaymath}
\end{theorem}
When $l=1$, we have
\begin{equation}
S_{1,\lambda}(n,1)\ =\ \sum_{k=1}^{n}S_{J,\lambda}^{(1)}(n,k)S_{2,\lambda}(k,1)\ =\ \sum_{k=1}^{n}(1)_{k,\lambda}S_{J,\lambda}^{(1)}(n,k).	\label{34}
\end{equation}
As an extension of the degenerate Bell polynomials in \eqref{12}, the Jindalrae polynomials are defined by
\begin{equation}
	e_{\lambda}^{x}\big(e_{\lambda}(e_{\lambda}(t)-1)-1\big)=\sum_{n=0}^{\infty}J_{n,\lambda}(x)\frac{t^{n}}{n!}\label{35}.
\end{equation}
When $x=1$, $J_{n,\lambda}=J_{n,\lambda}(1)$ are called Jindalrae numbers. \\
From \eqref{35}, we note that
\begin{align}
	e_{\lambda}^{x}\big(e_{\lambda}(e_{\lambda}(t)-1)-1\big) \ &=\ \sum_{k=0}^{\infty}(x)_{k,\lambda}\frac{1}{k!}\big(e_{\lambda}(e_{\lambda}(t)-1)-1\big)^{k} \label{36} \\
	&=\ \sum_{k=0}^{\infty}(x)_{k,\lambda}\sum_{n=k}^{\infty}S_{J,\lambda}^{(2)}(n,k)\frac{t^{n}}{n!} \nonumber \\
	&=\ \sum_{n=0}^{\infty}\bigg(\sum_{k=0}^{n}(x)_{k,\lambda}S_{J,\lambda}^{(2)}(n,k)\bigg)\frac{t^{n}}{n!}. \nonumber
\end{align}
Therefore, by \eqref{35} and \eqref{36}, we obtain the following theorem.
\begin{theorem}
	For $n\ge 0$, we have
	\begin{displaymath}
	J_{n,\lambda}(x)=\sum_{k=0}^{n}(x)_{k,\lambda}S_{J,\lambda}^{(2)}(n,k).
	\end{displaymath}
	In particular, for $x=1$,
	\begin{displaymath}
	J_{n,\lambda}=\sum_{k=0}^{n}(1)_{k,\lambda}S_{J,\lambda}^{(2)}(n,k).
	\end{displaymath}
\end{theorem}
By replacing $t$ by $\log_{\lambda}(1+t)$ in \eqref{35}, we get
\begin{align}
	e_{\lambda}^{x}\big(e_{\lambda}(t)-1\big)\ &=\ \sum_{m=0}^{\infty}J_{m,\lambda}(x)\frac{1}{m!}\big(\log_{\lambda}(1+t)\big)^{m} \label{37} \\
	&=\ \sum_{m=0}^{\infty}J_{m,\lambda}(x)\sum_{n=m}^{\infty}S_{1,\lambda}(n,m)\frac{t^{n}}{n!}\nonumber\\
	&=\ \sum_{n=0}^{\infty}\bigg(\sum_{m=0}^{n}J_{m,\lambda}(x)S_{1,\lambda}(n,m)\bigg)\frac{t^{n}}{n!}. \nonumber
\end{align}
Therefore, by \eqref{12} and \eqref{37}, we obtain the following theorem.
\begin{theorem}
	For $n\ge 0$, we have
	\begin{displaymath}
	B_{n,\lambda}(x)\ =\ \sum_{m=0}^{n}J_{m,\lambda}(x)S_{1,\lambda}(n,m).
	\end{displaymath}
	In particular,
	\begin{displaymath}
	B_{n,\lambda}\ =\ \sum_{m=0}^{n}J_{m,\lambda}S_{1,\lambda}(n,m).
	\end{displaymath}
\end{theorem}
From \eqref{20}, we note that
\begin{align}
	e_{\lambda}^{x}\big(e_{\lambda}(e_{\lambda}(t)-1)-1\big)\ &=\ \sum_{m=0}^{\infty}B_{m,\lambda}(x)\frac{1}{m!}\big(e_{\lambda}(t)-1\big)^{m} \label{38} \\
	&=\ \sum_{m=0}^{\infty}B_{m,\lambda}(x)\sum_{n=m}^{\infty}S_{2,\lambda}(n,m)\frac{t^{n}}{n!} \nonumber \\
	&=\ \sum_{n=0}^{\infty}\bigg(\sum_{m=0}^{n}B_{m,\lambda}(x)S_{2,\lambda}(n,m)\bigg)\frac{t^{n}}{n!}. \nonumber
\end{align}
Therefore, by \eqref{35} and \eqref{38}, we obtain the following theorem.
\begin{theorem}
	For $n\ge 0$, we have
	\begin{displaymath}
		J_{n,\lambda}(x)\ =\ \sum_{m=0}^{n}B_{m,\lambda}(x)S_{2,\lambda}(n,m).
	\end{displaymath}
	In particular,
	\begin{displaymath}
		J_{n,\lambda}\ =\ \sum_{m=0}^{n}B_{m,\lambda}S_{2,\lambda}(n,m).
	\end{displaymath}	
\end{theorem}
It is not difficult to show that $\log_{\lambda}\big(\log_{\lambda}(1+t)+1\big)$ is the compositional inverse of $e_{\lambda}\big(e_{\lambda}(t)-1\big)-1$. Now, we consider the Gaenari polynomials given by
\begin{equation}
	e_{\lambda}^{x}\big(\log_{\lambda}(\log_{\lambda}(1+t)+1)\big)\ =\ \sum_{n=0}^{\infty}G_{n,\lambda}(x)\frac{t^{n}}{n!}. \label{39}
\end{equation}
When $x=1$, $G_{n,\lambda}(1)=G_{n,\lambda}$ are called the Gaenari numbers. \\
From \eqref{6}, we note that
\begin{align}
	e_{\lambda}^{x}\big(\log_{\lambda}(\log_{\lambda}(1+t)+1)\big)\ &=\ \sum_{k=0}^{\infty}(x)_{k,\lambda}\frac{1}{k!} \big(\log_{\lambda}(\log_{\lambda}(1+t)+1)\big)^{k}\label{40} \\
	&=\ \sum_{k=0}^{\infty}(x)_{k,\lambda}\sum_{n=k}^{\infty}S_{J,\lambda}^{(1)}(n,k)\frac{t^{n}}{n!} \nonumber \\
	&=\ \sum_{n=0}^{\infty}\bigg(\sum_{k=0}^{n}(x)_{k,\lambda}S_{J,\lambda}^{(1)}(n,k)\bigg)\frac{t^{n}}{n!}. \nonumber
\end{align}
Therefore, by \eqref{39} and \eqref{40}, we obtain the following theorem.
\begin{theorem}
	For $n\ge 0$, we have
	\begin{displaymath}
		G_{n,\lambda}(x)\ =\ \sum_{k=0}^{n}(x)_{k,\lambda}S_{J,\lambda}^{(1)}(n,k).
	\end{displaymath}
	In particular,
	\begin{displaymath}
		G_{n,\lambda}\ =\ \sum_{k=0}^{n}(1)_{k,\lambda}S_{J,\lambda}^{(1)}(n,k).
	\end{displaymath}
\end{theorem}
By replacing $t$ by $e_{\lambda}(t)-1$ in \eqref{39}, we get
\begin{align}
	e_{\lambda}^{x}\big(\log_{\lambda}(1+t)\big)\ &= \sum_{m=0}^{\infty}G_{m,\lambda}(x)\frac{1}{m!}\big(e_{\lambda}(t)-1\big)^{m}\nonumber \\
	&=\ \sum_{m=0}^{\infty}G_{m,\lambda}(x)\sum_{n=m}^{\infty}S_{2,\lambda}(n,m)\frac{t^{n}}{n!} \label{41} \\
	&=\ \sum_{n=0}^{\infty}\bigg(\sum_{m=0}^{n}G_{m,\lambda}(x)S_{2,\lambda}(n,m)\bigg)\frac{t^{n}}{n!}. \nonumber
\end{align}
On the other hand,
\begin{equation}
	e_{\lambda}^{x}\big(\log_{\lambda}(1+t)\big)\ =\ \sum_{n=0}^{\infty}(x)_{n}\frac{t^{n}}{n!}. \label{42}
\end{equation}
Therefore, by \eqref{41} and \eqref{42}, we obtain the following theorem.
\begin{theorem}
	For $n\ge 0$, we have
	\begin{displaymath}
		(x)_{n}=\sum_{m=0}^{n}G_{m,\lambda}(x)S_{2,\lambda}(n,m).
	\end{displaymath}
	When $x=1$, we have
	\begin{equation}
		(1)_{n}=\sum_{m=0}^{n}G_{m,\lambda}S_{2,\lambda}(n,m).\label{43}
	\end{equation}
\end{theorem}
By \eqref{43}, we get
\begin{equation}
	\label{44} G_{0,\lambda}=1,\quad \sum_{m=0}^{n}G_{m,\lambda}S_{2,\lambda}(n,m)=\left\{\begin{array}{ccc}
	1, & \textrm{if $n=1$,}\\
	0, & \textrm{if $n > 1$.}
\end{array}\right.
\end{equation}
From \eqref{44}, we note that $G_{1,\lambda}=1$. Indeed, we note that
\begin{equation}
	e_{\lambda}^{x}\big(\log_{\lambda}(\log_{\lambda}(1+t)+1)\big)=\big(1+\log(1+t)\big)^{x}. \label{45}
\end{equation}
Thus, by \eqref{39} and \eqref{45}, we get
\begin{equation}
	\big(1+\log_{\lambda}(1+t)\big)^{x}=\sum_{n=0}^{\infty}G_{n,\lambda}(x)\frac{t^{n}}{n!}. \label{46}
\end{equation}
From \eqref{7} and \eqref{46}, we have
\begin{equation}
	\sum_{n=0}^{\infty}G_{n,\lambda}\frac{t^{n}}{n!}=1+\sum_{n=1}^{\infty}\lambda^{n-1}(1)_{n,1/\lambda}\frac{t^{n}}{n!}. \label{47}
\end{equation}
Thus, by \eqref{47}, we get
\begin{displaymath}
	G_{n,\lambda}=\lambda^{n-1}(1)_{n,1/\lambda},\quad (n\ge 1).
\end{displaymath}
\begin{corollary}
	For $n\ge 1$, we have
	\begin{displaymath}
			G_{n,\lambda}=\lambda^{n-1}(1)_{n,1/\lambda},\quad (n\ge 1).
	\end{displaymath}
\end{corollary}
By replacing $t$ by $e_{\lambda}\big(e_{\lambda}(t)-1\big)-1$ in \eqref{39}, we get
\begin{align}
	e_{\lambda}^{x}(t)\ &=\ \sum_{m=0}^{\infty}G_{m,\lambda}(x)\frac{1}{m!}\big(e_{\lambda}(e_{\lambda}(t)-1\big)^{m}\label{48}\\
	&=\ \sum_{m=0}^{\infty}G_{m,\lambda}(x)\sum_{n=m}^{\infty}S_{J,\lambda}^{(2)}(n,m)\frac{t^{n}}{n!}\nonumber\\
	&=\ \sum_{n=0}^{\infty}\bigg(\sum_{m=0}^{n}G_{m,\lambda}(x)S_{J,\lambda}^{(2)}(n,m)\bigg)\frac{t^{n}}{n!}. \nonumber
\end{align}
From \eqref{6} and \eqref{48}, we have
\begin{equation}
	(x)_{n,\lambda}=\sum_{m=0}^{n}G_{m,\lambda}(x)S_{J,\lambda}^{(2)}(n,m),\quad (n\ge 0). \label{49}
\end{equation}
In particular,
\begin{displaymath}
	(1)_{n,\lambda}=\sum_{m=0}^{n}G_{m,\lambda}S_{J,\lambda}^{(2)}(n,m).
\end{displaymath}
From \eqref{35}, we also note that
\begin{align}
	e_{\lambda}^{x}(t)\ &=\ \sum_{m=0}^{\infty}J_{m,\lambda}(x)\frac{1}{m!}\big(\log_{\lambda}(\log_{\lambda}(1+t)+1)\big)^{m} \label{50} \\
	&=\ \sum_{m=0}^{\infty}J_{m,\lambda}(x)\sum_{n=m}^{\infty}S_{J,\lambda}^{(1)}(n,m)\frac{t^{n}}{n!} \nonumber \\
	&=\ \sum_{n=0}^{\infty}\bigg(\sum_{m=0}^{n}J_{m,\lambda}(x)S_{J,\lambda}^{(1)}(n,m)\bigg)\frac{t^{n}}{n!}\nonumber.
\end{align}
Thus, by \eqref{6} and \eqref{50}, we get
\begin{equation}
	(x)_{n,\lambda}=\sum_{m=0}^{n}J_{m,\lambda}(x)S_{J,\lambda}^{(1)}(n,m),\quad(n\ge 0). \label{51}
\end{equation}
From \eqref{49} and \eqref{51}, we have
\begin{equation}
	\sum_{m=0}^{n}G_{m,\lambda}(x)S_{J,\lambda}^{(2)}(n,m)\ =\ \sum_{m=0}^{n}J_{m,\lambda}(x)S_{J,\lambda}^{(1)}(n,m).\label{52}
\end{equation}

\vspace{0.1 in}

\section{Further remarks}
In this section, we are going to interpret what we obtained in the previous section by means of umbral calculus.
First, we will go over some necessary facts about umbral calculus. For more details on umbral calculus, we let the reader refer to [28]. \\
\indent A series $g(t)$ with $O(g(t))=0$ and a series $f(t)$ with $O(f(t))=1$ are respectively called an invertible series and a delta series. Recall here that the order $O(f(t))$ of the non-zero power series $f(t)$ is the smallest integer $k$ for which the coefficient of $t^k$ does not vanish.
Let $g(t)$ be an invertible series and let $f(t)$ be a delta series. Then there exists a unique sequence $s_n(x)$ of polynomials such that $\langle g(t)f(t)^k|s_n(x)\rangle=n!\delta_{n,k}$, for $n,k\geq0$ (see [28]). The sequence $s_n(x)$ is called the Sheffer sequence for the Sheffer pair $(g(t),f(t))$, which is denoted by $s_n(x)\sim(g(t),f(t))$. In particular, $s_n(x)$ is called the associated sequence to $f(t)$, if $s_n(x)\sim(1,f(t))$. Further,  $s_n(x)$ is called the Appell sequence for $g(t)$, if $s_n(x)\sim(g(t),t)$.  It is well known that $s_n(x)\sim(g(t),f(t))$ if and only if
\begin{equation}\label{53}
\frac{1}{g(\bar{f}(t))}e^{x\bar{f}(t)}=\sum_{n\geq0}s_n(x)\frac{t^n}{n!},
\end{equation}
where $\bar{f}(t)$ is the compositional inverse of $f(t)$ determined by $f(\bar{f}(t))=\bar{f}(f(t))=t$.

For each nonnegative integer $m$, the $m$th power of an invertible series $g(t)$ will be indicated by $(g(t))^m$, while the compositional powers of a delta series $f(t)$ will be denoted by $f^m(t)=f\circ f\circ\cdots\circ f(t)$\,\,($m$ times). For $p_n(x)$ and $q_n(x)=\sum_{k=0}^nq_{n,k}t^k$, the umbral composition of $q_n(x)$ with $p_n(x)$, denoted by $q_n\circ p_n(x)$, is defined by $q_n\circ p_n(x)=\sum_{k=0}^n q_{n,k}p_k(x)$.\\
\indent The next result is stated in Theorem 3.5.5 of [28].

\begin{theorem}
The set of Sheffer sequences forms a group under operation of umbral composition. If $s_n(x)\sim(g(t),f(t))$ and $r_n(x)\sim(h(t),\ell(t))$, then $r_n(x)\circ s_n(x)\sim(g(t)h(f(t)),\ell(f(t)))$. The identity under umbral composition is $x^n\sim(1,t)$, and the inverse of the sequence $s_n(x)\sim(g(t),f(t))$ is the Sheffer sequence for $(g(\bar{f}(t))^{-1},\bar{f}(t))$.
\end{theorem}

As a corollary, we obtain the following result that will be needed later.

\begin{corollary}
Let $s_n(x)\sim(g(t),f(t))$, and let $r_n(x)\sim(1,\ell(t))$. Then, for any positive integer $m$, the generating function for $r_n^{(m)}\circ s_n(x)$ is obtained from that for $s_n(x)$ by substituting $\bar{\ell}^m(t)$ for $t$.
\begin{proof}
As $r_n^{(m)}\circ s_n(x)\sim(g(t),\ell^m(f(t)))$, and the compositional inverse of $\ell^m(f(t))$ is $\bar{f}(\bar{\ell}^m(t))$, we have
\begin{align}\label{54}
g(\bar{f}(\bar{\ell}^m(t)))^{-1}e^{x\bar{f}(\bar{\ell}^m(t))}=\sum_{n\geq0}r_n^{(m)}\circ s_n(x)\frac{t^n}{n!}.
\end{align}
\end{proof}
\end{corollary}

From the definition of umbral composition, we see that the $m$th power under umbral composition of $r_n(x)\sim(h(t),\ell(t))$ is given by
\begin{align}\label{55}
r_n^{(m)}(x)\sim\left(\prod_{i=1}^{m-1}h(\ell^i(t)),\ell^m(t)\right),
\end{align}
for any positive integer $m$. In particular, this says that, for the Appell sequence $r_n(x)\sim(h(t),t)$, we have $r_n^{(m)}(x)\sim((h(t))^m,t)$. Whereas, for the associated sequence $r_n(x)\sim(1,\ell(t))$, we have $r_n^{(m)}(x)\sim(1,\ell^m(t))$.

For $n\geq0$, we write $r_n(x)=\sum_{k=0}^n r_{n,k}x^k=\sum_{k\geq0}r_{n,k}x^k$, where we agree that $r_{i,j}=0$ for all $i<j$. In general, we write
$$r_n^{(m)}(x)=\sum_{k-0}^nr_{n,k}^{(m)}x^k=\sum_{k\geq0}r_{n,k}^{(m)}x^k,$$
for all $m\in\mathbb{Z}_{>0}$. Then we see that
\begin{align}\label{56}
r^{(m)}_{n,k}&=\sum_{\ell_1,\ldots,\ell_{m-1}=0}^nr_{n,\ell_1}r_{\ell_1,\ell_2}\cdots r_{\ell_{m-1},k},\quad m\geq1,
\end{align}
where we understand $r_{n,k}^{(1)}=r_{n,k}$, for $m=1$.

\subsection{Jindalrae polynomials}
From \eqref{35}, we note that the Jindalrae polynomials $J_{n,\lambda}(x)$ are given by
\begin{equation}\label{57}
\sum_{n=0}^{\infty}J_{n,\lambda}(x)\frac{t^n}{n!}=e^{x\bar{f}(\bar{\ell}^2(t))},
\end{equation}
where $f(t)=\frac{1}{\lambda}(e^{\lambda t} -1), \,\,\ell(t)=\log_{\lambda}(1+t)$. \\
\indent Noting that $\bar{f}(t)=\frac{1}{\lambda}\log(1+\lambda t), \bar{\ell}(t)=e_{\lambda}(t)-1$, we have
\begin{equation}\begin{split}\label{58}
&e^{x\bar{f}(t)}=e_{\lambda}^{x}(t)=\sum_{n=0}^{\infty}(x)_{n,\lambda}\frac{t^n}{n!}, \\
&e^{x\bar{\ell}(t)}=\sum_{n=0}^{\infty}\sum_{k=0}^{n}S_{2,\lambda}(n,k)x^k\frac{t^n}{n!}.
\end{split}\end{equation}
By \eqref{58}, we let
\begin{equation}\begin{split}\label{59}
&s_n(x)=(x)_{n,\lambda}\sim (1,f(t)=\frac{1}{\lambda}(e^{\lambda t}-1)), \\
&r_n(x)=\sum_{k=0}^{n}S_{2,\lambda}(n,k)x^k\sim(1,\ell(t)=\log_{\lambda}(1+t)).
\end{split}\end{equation}
From \eqref{54} and \eqref{57}, we observe that
\begin{equation}\label{60}
J_{n,\lambda}(x)=r_n^{(2)}\circ s_n(x).
\end{equation}
As a check, we verify Theorem 8 again by using \eqref{60}. From \eqref{56} and \eqref{59}, we see that
\begin{equation}\begin{split}\label{61}
r_{n}^{(2)}(x)&=\sum_{k=0}^{n}\sum_{l=0}^{n}S_{2,\lambda}(n,l)S_{2,\lambda}(l,k)x^k \\
&=\sum_{k=0}^{n}\sum_{l=k}^{n}S_{2,\lambda}(n,l)S_{2,\lambda}(l,k)x^k \\
&=\sum_{k=0}^{\infty}S_{J,\lambda}^{(2)}(n,k)x^k.
\end{split}\end{equation}
Thus, from \eqref{59}, \eqref{60} and \eqref{61},  we get the desired result as follows:
$$ J_{n, \lambda}(x)=\sum_{k=0}^{n}S_{J, \lambda}^{(2)}(n,k)(x)_{k,\lambda}.$$
\indent  For any positive integer $r$, Korobov polynomials of the first kind of order $r$ are given by
\begin{equation}\label{62}
\bigg(\frac{t}{\log_{\lambda}(1+t)}\bigg)^r(1+t)^x=\bigg(\frac{\lambda t}{(1+t)^{\lambda}-1}\bigg)^r(1+t)^x=\sum_{n=0}^{\infty}K_{n, (r)}(x;\lambda).
\end{equation}
For $x=0$, $K_{n, (r)}(\lambda)=K_{n, (r)}(0;\lambda)$ are called Korobov numbers of the first kind of order $r$. The Korobov polynomials (respectively, numbers) of the first kind are also called the degenerate Bernoulli polynomials (respectively, numbers) of the second. \\
\indent In [9], it was shown that
\begin{equation*}\begin{split}
&\sum_{l_1,\dots,l_{m-1}=0}^{n}S_{2,\lambda}(n,l_1)S_{2,\lambda}(l_1,l_2) \cdots S_{2,\lambda}(l_{m-1},k) \\
&=\sum_{n\geq l_1 \geq \dots \geq l_{m-1}\geq k}S_{2,\lambda}(n,l_1)S_{2,\lambda}(l_1,l_2) \cdots S_{2,\lambda}(l_{m-1},k) \\
&=\sum_{k_1 +\cdots+ k_m=n-k}\binom{n-1}{k_1,\dots, k_m, k-1}\prod_{j=1}^mK_{k_j,(n-\sum_{i=j+1}^mk_i)}(\lambda).
\end{split}\end{equation*}
Therefore, in the special cases of $m=1$ and $m=2$, we respectively have
\begin{equation*}\begin{split}
&S_{2,\lambda}(n,k)=\binom{n-1}{k-1}K_{n-k, (n)}(\lambda), \\
&S_{J,\lambda}^{(2)}(n,k)=\sum_{k_1+k_2=n-k}\binom{n-1}{k_1,k_2,k-1}K_{k_2,(n)}(\lambda)K_{k_1,(n-k_2)}(\lambda).
\end{split}\end{equation*}

\subsection{Gaenari polynomials}

From \eqref{39}, we recall that the Gaenari polynomials $G_{n,\lambda}(x)$ are given by
\begin{equation}\label{63}
\sum_{n=0}^{\infty}G_{n,\lambda}(x)\frac{t^n}{n!}=e^{x\bar{f}(\bar{\ell}^2(t))},
\end{equation}
where $f(t)=\frac{1}{\lambda}(e^{\lambda t} -1), \,\,\ell(t)=e_{\lambda}(t)-1$. \\
\indent Noting that $\bar{f}(t)=\frac{1}{\lambda}\log(1+\lambda t), \bar{\ell}(t)=\log_{\lambda}(1+t)$, we have
\begin{equation}\begin{split}\label{64}
&e^{x\bar{f}(t)}=e_{\lambda}^{x}(t)=\sum_{n=0}^{\infty}(x)_{n,\lambda}\frac{t^n}{n!}, \\
&e^{x\bar{\ell}(t)}=\sum_{n=0}^{\infty}\sum_{k=0}^{n}S_{1,\lambda}(n,k)x^k\frac{t^n}{n!}.
\end{split}\end{equation}
By \eqref{58}, we let
\begin{equation}\begin{split}\label{65}
&s_n(x)=(x)_{n,\lambda}\sim (1,f(t)=\frac{1}{\lambda}(e^{\lambda t}-1)), \\
&r_n(x)=\sum_{k=0}^{n}S_{1,\lambda}(n,k)x^k\sim(1,\ell(t)=e_{\lambda}(t)-1).
\end{split}\end{equation}
From \eqref{54} and \eqref{63}, we observe that
\begin{equation}\label{66}
G_{n,\lambda}(x)=r_n^{(2)}\circ s_n(x).
\end{equation}
Here we would like to verify Theorem 11 again by making use of \eqref{66}. From \eqref{56} and \eqref{65}, we see that
\begin{equation}\begin{split}\label{67}
r_{n}^{(2)}(x)&=\sum_{k=0}^{n}\sum_{l=0}^{n}S_{1,\lambda}(n,l)S_{1,\lambda}(l,k)x^k \\
&=\sum_{k=0}^{n}\sum_{l=k}^{n}S_{1,\lambda}(n,l)S_{1,\lambda}(l,k)x^k \\
&=\sum_{k=0}^{\infty}S_{J,\lambda}^{(1)}(n,k)x^k.
\end{split}\end{equation}
Thus, from \eqref{65}, \eqref{66} and \eqref{67},  we get what we wanted as follows:
$$ G_{n, \lambda}(x)=\sum_{k=0}^{n}S_{J, \lambda}^{(1)}(n,k)(x)_{k,\lambda}.$$
\indent  For any positive integer $r$, the degenerate Bernoulli  polynomials $\beta_{n,(r)}(x;\lambda)$ of order $r$ are defined by
\begin{equation}\label{68}
\bigg(\frac{t}{e_{\lambda}(t)-1}\bigg)^re_{\lambda}^x(t)=\sum_{n\geq0}\beta_{n,(r)}(x;\lambda)\frac{t^n}{n!}.
\end{equation}
For $x=0$, $\beta_{n, (r)}(\lambda)=\beta_{n, (r)}(0;\lambda)$ are called the degenerate Bernoulli numbers of order $r$.  \\
\indent In [9], it was shown that
\begin{equation*}\begin{split}
&\sum_{l_1,\dots,l_{m-1}=0}^{n}S_{1,\lambda}(n,l_1)S_{1,\lambda}(l_1,l_2) \cdots S_{1,\lambda}(l_{m-1},k) \\
&=\sum_{n\geq l_1 \geq \dots \geq l_{m-1}\geq k}S_{1,\lambda}(n,l_1)S_{1,\lambda}(l_1,l_2) \cdots S_{1,\lambda}(l_{m-1},k) \\
&=\sum_{k_1+\cdots+k_m=n-k}\binom{n-1}{k_1,\ldots,k_m,k-1}\prod_{j=1}^m\beta_{k_j,(n-\sum_{i=j+1}^mk_i)}(\lambda).
\end{split}\end{equation*}
Hence, in the special cases of $m=1$ and $m=2$, we respectively have
\begin{equation*}\begin{split}
&S_{1,\lambda}(n,k)=\binom{n-1}{k-1}\beta_{n-k, (n)}(\lambda), \\
&S_{J,\lambda}^{(1)}(n,k)=\sum_{k_1+k_2=n-k}\binom{n-1}{k_1,k_2,k-1}\beta_{k_2,(n)}(\lambda)\beta_{k_1,(n-k_2)}(\lambda).
\end{split}\end{equation*}

\vspace{0.1in}

\section{Conclusion}
In this paper, we introduced Jindalrae-Stirling numbers of the first and second kinds, as extensions of the notions of the degenerate Stirling numbers of the first and second kinds and found some relations connecting those special numbers, the degenerate Stirling numbers of the first and second kinds and the degenerate Bell numbers and polynomials. Then we introduced Jindalrae numbers and polynomials, as an extension of the notion of the degenerate Bell numbers and polynomials, and Gaenari numbers and polynomials, and obtained some explicit expressions and identities involving those numbers and polynomials, Jindalrae-Stirling numbers of the first and second kinds, the degenerate Stirling numbers of the first and second kinds and the degenerate Bell polynomials.  In addition, we interpreted our results by means of umbral calculus. \\
\indent As to possible applications our results, we would like to mention three things. The first one is their applications to differential equations. In [8], new combinatorial identities for some degenerate special polynomials were found from certain infinite families of linear and non-linear ordinary differential equations, satisfied by the generating functions of those polynomials. The second one is their applications to probability theory. In [16,19], new identities connecting some special numbers and moments
of random variables were derived by using the generating functions of the moments of certain random variables.
The last one is their applications to identities of symmetry. In [14], abundant identities of symmetry were derived  for various degenerate versions of many special polynomials by using $p$-adic integrals. \\
\indent It is one of our future projects to continue this line of research, namely study of degenerate versions of some special polynomials and numbers, and to find some of their possible applications to mathematics, science and engineering.

\end{document}